\documentclass[11pt,reqno]{amsart}
\usepackage{amsmath,amssymb,geometry,color,tikz-cd}
\usepackage[backref,pagebackref,pdftex,hyperindex]{hyperref}
\usepackage[alphabetic]{amsrefs}
\usepackage{amssymb}% see geometry.pdf on how to lay out the page. There's lots.
\usepackage{bbm}
\usepackage{enumitem}
\usepackage{upgreek}

\newtheorem{proposition}{Proposition}[section]
\newtheorem{theorem}[proposition]{Theorem}
\newtheorem{lemma}[proposition]{Lemma}

\newtheorem{conjecture}[proposition]{Conjecture}

\theoremstyle{definition}
\newtheorem{remark}[proposition]{Remark}
\newtheorem{definition}[proposition]{Definition}
\newtheorem{example}[proposition]{Example}

\title{On K-semistable domains - more examples}
\author{Chuyu Zhou}
\address{Department of Mathematics, Yonsei University, 50 Yonsei-ro, Seodaemun-gu, Seoul 03722, Republic of Korea}
\email{chuyuzhou1@gmail.com}
\date{} % delete this line to display the current date

\thanks{
        \newline
		\indent
		\emph{Competing interests}: The author declares none.
		}

\newcommand{\ord}{{\rm {ord}}}
\newcommand{\tc}{{\rm {tc}}}
\newcommand{\vol}{{\rm {vol}}}

\newcommand{\dt}{{\rm {dt}}}

\newcommand{\Kss}{{\rm {Kss}}}

\newcommand{\bA}{\mathbb{A}}

\newcommand{\bC}{\mathbb{C}}

\newcommand{\bN}{\mathbb{N}}

\newcommand{\bP}{\mathbb{P}}
\newcommand{\bQ}{\mathbb{Q}}

\newcommand{\mE}{\mathcal{E}}

\newcommand{\mO}{\mathcal{O}}

\newcommand{\mV}{\mathcal{V}}

\newcommand{\mX}{\mathcal{X}}

\begin{document}

\begin{abstract}
We compute K-semistable domains for various examples of log pairs.
\end{abstract}

\maketitle
\tableofcontents

\section{Introduction}

We work over the complex number field $\bC$ throughout the article. 

In \cite{LZ23}, K-semistable domains (see Definition \ref{def: kss domain} and Remark \ref{rem: kss domain}) for various examples were computed. In this paper, we apply the approach in \cite{LZ23} to compute K-semistable domains for more general examples, which recover all the examples computed there. We list these examples as the following five theorems. In fact, Theorem \ref{thm: main PSS} is a special case of Theorem \ref{thm: main SSS}, Theorem \ref{thm: main SSS} is a special case of Theorem \ref{thm: main SSS|S}, and Theorem \ref{thm: main PSSS} is a special case of Theorem \ref{thm: main SSSS}. However, for the readability and the convenience of checking the computation, we prove all of them, and in fact we sometimes need to use the special case to prove the general case.

\textit{Notation}: Let $(\bP^n, \sum_{j=1}^k S_{d_j})$ be a log smooth pair such that $S_{d_j}, j=1,...,k,$ are mutually different smooth hypersurfaces of degrees $d_j$. We say all the Fano complete intersections are K-semistable if for any non-repeating subset $\{i_{j_1},...,i_{j_l}\}\subset \{1, 2,...,k\}$ satisfying that $\cap_{s=1}^l S_{d_{j_s}}$ is a Fano manifold, the complete intersection $\cap_{s=1}^l S_{d_{j_s}}$
is K-semistable.

Before stating the results, we first recall the concept of K-semistable domain.
Let $X$ be a Fano variety and $D_j, 1\leq j\leq k,$ are effective $\bQ$-divisors on $X$ such that each $D_j$ is proportional to $-K_X$. Then we define the K-semistable domain for $(X, \sum_{j=1}^kD_j)$ as follows:
$$\Kss(X, \sum_{j=1}^kD_j):=\overline{\{(x_1,...,x_k)\in \bQ_{\geq 0}^k\ |\ \textit{$(X, \sum_{j=1}^kx_jD_j)$ is a K-semistable log Fano pair}\}},$$
where the overline means taking the closure.

\begin{theorem}\label{thm: main PSS}{\rm{(Theorem \ref{thm: PSS}, Example \ref{exa: PSS})}}
Consider the log pair $(\bP^n, S_{d_1}+S_{d_2})$, where $S_{d_1}, S_{d_2}$ are two general smooth hypersurfaces of degrees $d_1, d_2\leq n+1$ and $n\geq 2$. Suppose $S_{d_1}\cap S_{d_2}$ is K-semistable if $d_1+d_2\leq n+1$.
\iffalse
Suppose the following two conditions:
\begin{enumerate}
\item $(\bP^n, S_{d_1}+S_{d_2})$ is log canonical,
\item $S_{d_1}\cap S_{d_2}$ is K-semistable of dimension $n-2$ if $d_1+d_2\leq n+1$. 
\end{enumerate}
\fi
Then 
$\Kss(\bP^n, S_{d_1}+S_{d_2})$ is a polytope $P$ generated by the following equations
\begin{equation*}
\begin{cases}
0\leq x\leq 1\\
0\leq y\leq 1\\
d_1x+d_2y\leq n+1\\
\frac{d_1-1}{d_1}-\frac{nx}{n+1}+\frac{d_2y}{d_1(n+1)}\geq 0\\
\frac{d_2-1}{d_2}-\frac{ny}{n+1}+\frac{d_1x}{d_2(n+1)}\geq 0.\\
\end{cases}
\end{equation*}
If $d_1+d_2<n+1$, then $P$ is generated by the extremal points 
$$(0,0), \ \ \left(\frac{(n+1)(d_1-1)}{nd_1}, 0\right),\ \ \left(0, \frac{(n+1)(d_2-1)}{nd_2}\right), \ \ \left(\frac{d_2-1+n(d_1-1)}{(n-1)d_1}, \frac{d_1-1+n(d_2-1)}{(n-1)d_2}\right). $$
If $d_1+d_2\geq n+1$, then $P$ is generated by the extremal points
$$(0,0), \ \ \left(\frac{(n+1)(d_1-1)}{nd_1}, 0\right),\ \ \left(0, \frac{(n+1)(d_2-1)}{nd_2}\right), \ \ \left(1, \frac{n+1-d_1}{d_2}\right), \ \ \left(\frac{n+1-d_2}{d_1}, 1\right). $$ 
\end{theorem}

\begin{theorem}\label{thm: main SSS}{\rm{(Theorem \ref{thm: SSS}, Example \ref{exa: SSS})}}
Consider the log pair $(S_d, {S_{d_1}}|_{S_d}+{S_{d_2}}|_{S_d})$, where $S_d, S_{d_1}, S_{d_2}$ are three general smooth hypersurfaces in $\bP^{n+1}$ of degrees $ d, d_1, d_2$ with $n\geq 2$ and $ d_1, d_2\leq n+2-d$. 
Suppose all the Fano complete intersections are K-semistable.
%Suppose $S_d\cap S_{d_1}\cap S_{d_2}$ is K-semistable of dimension $n-2$ if $d_1+d_2\leq n+2-d$.
\iffalse
Suppose the following two conditions:
\begin{enumerate}
\item $(S_d, {S_{d_1}}|_{S_d}+{S_{d_2}}|_{S_d})$ is log canonical,
\item $S_d\cap S_{d_1}\cap S_{d_2}$ is K-semistable of dimension $n-2$ if $d_1+d_2\leq n+2-d$. 
\end{enumerate}
\fi
Then $\Kss(S_d, {S_{d_1}}|_{S_d}+{S_{d_2}}|_{S_d})$ is a polytope $P$ generated by the following equations
\begin{equation*}
\begin{cases}
0\leq x\leq 1\\
0\leq y\leq 1\\
d_1x+d_2y\leq n+2-d\\
\frac{d_1(n+1)-(n+2-d)}{d_1(n+1)}-\frac{nx}{n+1}+\frac{d_2y}{d_1(n+1)}\geq 0\\
\frac{d_2(n+1)-(n+2-d)}{d_2(n+1)}-\frac{ny}{n+1}+\frac{d_1x}{d_2(n+1)}\geq 0.\\
\end{cases}
\end{equation*}
If $d_1+d_2<n+2-d$, then $P$ is generated by the extremal points 
$$(0,0), \ \ \left(\frac{(n+1)d_1-(n+2-d)}{nd_1}, 0\right),\ \ \left(0, \frac{(n+1)d_2-(n+2-d)}{nd_2}\right)$$
and
$$\left(\frac{d_2+nd_1-(n+2-d)}{(n-1)d_1}, \frac{d_1+nd_2-(n+2-d)}{(n-1)d_2}\right). $$
If $d_1+d_2\geq n+2-d$, then $P$ is generated by the extremal points
$$(0,0), \ \ \left(\frac{(n+1)d_1-(n+2-d)}{nd_1}, 0\right),\ \ \left(0, \frac{(n+1)d_2-(n+2-d)}{nd_2}\right)$$
and
$$\left(1, \frac{n+2-d-d_1}{d_2}\right), \ \ \left(\frac{n+2-d-d_2}{d_1}, 1\right). $$ 
\end{theorem}

\begin{theorem}\label{thm: main SSS|S}{\rm{(Theorem \ref{thm: SSS|S}, Example \ref{exa: SSS|S})}}
Consider the log pair 
$$\left({S_{d_1}}|_{S_d}, {({S_{d_2}}|_{S_d})}|_{{S_{d_1}}|_{S_d}}+{({S_{d_3}}|_{S_d})}|_{{S_{d_1}}|_{S_d}}\right),$$ 
where $S_d, S_{d_1}, S_{d_2}$ and $S_{d_3}$ are four general smooth hypersurfaces of degrees $ d, d_1, d_2, d_3$ in $\bP^{n+2}$ with $n\geq 2$ and $d_2, d_3\leq n+3-d-d_1$.  
Suppose all the Fano complete intersections are K-semistable.
%Suppose $S_d\cap S_{d_1}\cap S_{d_2}\cap S_{d_3}$ is K-semistable of dimension $n-2$ if $d_2+d_3\leq n+3-d-d_1$.
\iffalse
Suppose the following two conditions:
\begin{enumerate}
\item $\left({S_{d_1}}|_{S_d}, {({S_{d_2}}|_{S_d})}|_{{S_{d_1}}|_{S_d}}+{({S_{d_3}}|_{S_d})}|_{{S_{d_1}}|_{S_d}}\right)$ is log canonical,
\item $S_d\cap S_{d_1}\cap S_{d_2}\cap S_{d_3}$ is K-semistable of dimension $n-2$ if $d_2+d_3\leq n+3-d-d_1$. 
\end{enumerate}
\fi
Then 
$$\Kss\left({S_{d_1}}|_{S_d}, {({S_{d_2}}|_{S_d})}|_{{S_{d_1}}|_{S_d}}+{({S_{d_3}}|_{S_d})}|_{{S_{d_1}}|_{S_d}}\right)$$
is  a polytope $P$ generated by the following equations
\begin{equation*}
\begin{cases}
0\leq x\leq 1\\
0\leq y\leq 1\\
d_2x+d_3y\leq n+3-d-d_1\\
\frac{d_2(n+1)-(n+3-d-d_1)}{d_2(n+1)}-\frac{nx}{n+1}+\frac{d_3y}{d_2(n+1)}\geq 0\\
\frac{d_3(n+1)-(n+3-d-d_1)}{d_3(n+1)}-\frac{ny}{n+1}+\frac{d_2x}{d_3(n+1)}\geq 0.\\
\end{cases}
\end{equation*}
If $d_2+d_3<n+3-d-d_1$, then $P$ is generated by the extremal points 
$$(0,0), \ \ \left(\frac{(n+1)d_2-(n+3-d-d_1)}{nd_2}, 0\right),\ \ \left(0, \frac{(n+1)d_3-(n+3-d-d_1)}{nd_3}\right)$$
and
$$\left(\frac{nd_2+d_1+d_3-(n+3-d)}{(n-1)d_2}, \frac{nd_3+d_1+d_2-(n+3-d)}{(n-1)d_3}\right). $$
If $d_2+d_3\geq n+3-d-d_1$, then $P$ is generated by the extremal points
$$(0,0), \ \ \left(\frac{(n+1)d_2-(n+3-d-d_1)}{nd_2}, 0\right),\ \ \left(0, \frac{(n+1)d_3-(n+3-d-d_1)}{nd_3}\right)$$
and 
$$\left(1, \frac{n+3-d-d_1-d_2}{d_3}\right),\ \  \left(\frac{n+3-d-d_1-d_3}{d_2}, 1\right). $$
\end{theorem}

\begin{theorem}\label{thm: main PSSS}{\rm{(Theorem \ref{thm: PSSS}, Example \ref{exa: PSSS})}}
Consider the log pair $(\bP^n, S_{d_1}+S_{d_2}+S_{d_3})$, where $S_{d_1}, S_{d_2}, S_{d_3}$ are three general smooth hypersurfaces of degrees $d_1, d_2, d_3$ in $\bP^{n}$ with $n\geq 2$ and $d_1, d_2, d_3\leq n+1$.  
Suppose all the Fano complete intersections are K-semistable.
%Suppose $S_{d_1}\cap S_{d_2}\cap S_{d_3}$ is K-semistable of dimension $n-3$ if $d_1+d_2+d_3\leq n+1$. 
\iffalse
Suppose the following two conditions:
\begin{enumerate}
\item $(\bP^n, S_{d_1}+S_{d_2}+S_{d_3})$ is log canonical,
\item $S_{d_1}\cap S_{d_2}\cap S_{d_3}$ is K-semistable of dimension $n-3$ if $d_1+d_2+d_3\leq n+1$. 
\end{enumerate}
\fi
Then 
$\Kss(\bP^n, S_{d_1}+S_{d_2}+S_{d_3})$ is a polytope $P$ generated by the following equations
\begin{equation*}
\begin{cases}
0\leq x\leq 1\\
0\leq y\leq 1\\
0\leq z\leq 1\\
d_1x+d_2y+d_3z\leq n+1\\
\frac{d_1-1}{d_1}-\frac{nx}{n+1}+\frac{d_2y}{d_1(n+1)}+\frac{d_3z}{d_1(n+1)}\geq 0\\
\frac{d_2-1}{d_2}-\frac{ny}{n+1}+\frac{d_1x}{d_2(n+1)}+\frac{d_3z}{d_2(n+1)}\geq 0
\\
\frac{d_3-1}{d_3}-\frac{nz}{n+1}+\frac{d_1x}{d_3(n+1)}+\frac{d_2y}{d_3(n+1)}\geq 0.
\end{cases}
\end{equation*}
If $d_1+d_2+d_3<n+1$, then $P$ is generated by the extremal points 
$$(0,0,0), \ \ \left(\frac{(n+1)(d_1-1)}{nd_1}, 0,0\right),\ \ \left(0, \frac{(n+1)(d_2-1)}{nd_2}, 0\right), \ \ \left(0,0, \frac{(n+1)(d_3-1)}{nd_3}\right),$$

$$\left(0, \frac{d_3-1+n(d_2-1)}{(n-1)d_2}, \frac{d_2-1+n(d_3-1)}{(n-1)d_3}\right),$$

$$\left(\frac{d_3-1+n(d_1-1)}{(n-1)d_1}, 0,  \frac{d_1-1+n(d_3-1)}{(n-1)d_3}\right),$$

$$\left(\frac{d_2-1+n(d_1-1)}{(n-1)d_1}, \frac{d_1-1+n(d_2-1)}{(n-1)d_2}, 0\right),
$$

$$\left(\frac{d_2+d_3+(n-1)d_1-(n+1)}{(n-2)d_1}, \frac{d_1+d_3+(n-1)d_2-(n+1)}{(n-2)d_2}, \frac{d_1+d_2+(n-1)d_3-(n+1)}{(n-2)d_3}\right). $$
If $d_1+d_2+d_3\geq n+1$, then $P$ is generated by the extremal points
$$(0,0,0), \ \ \left(\frac{(n+1)(d_1-1)}{nd_1}, 0,0\right),\ \ \left(0, \frac{(n+1)(d_2-1)}{nd_2}, 0\right), \ \ \left(0,0, \frac{(n+1)(d_3-1)}{nd_3}\right),$$
and $A_1 ,B_1, A_2 ,B_2, A_3 ,B_3$, where
\begin{equation*}
A_1, B_1=
\begin{cases}
(0, 1, \frac{n+1-d_2}{d_3}),\ (0,\frac{n+1-d_3}{d_2}, 1) & \text{if \ $d_2+d_3\geq n+1$}\\
\ \\
\left(0, \frac{d_3-1+n(d_2-1)}{(n-1)d_2}, \frac{d_2-1+n(d_3-1)}{(n-1)d_3}\right),\ \left(\frac{n+1-(d_2+d_3)}{d_1}, 1,1\right) & \text{if \ $d_2+d_3\leq n+1$}\\
\end{cases}
\end{equation*}

\begin{equation*}
A_2, B_2=
\begin{cases}
(1, 0, \frac{n+1-d_1}{d_3}),\ (\frac{n+1-d_3}{d_1}, 0, 1) & \text{if \ $d_1+d_3\geq n+1$}\\
\ \\
\left(\frac{d_3-1+n(d_1-1)}{(n-1)d_1}, 0,  \frac{d_1-1+n(d_3-1)}{(n-1)d_3}\right),\ \left(1, \frac{n+1-(d_1+d_3)}{d_2}, 1\right) & \text{if \ $d_1+d_3\leq n+1$}\\
\end{cases}
\end{equation*}

\begin{equation*}
A_3, B_3=
\begin{cases}
(1, \frac{n+1-d_1}{d_2}, 0),\ (\frac{n+1-d_2}{d_1}, 1, 0) & \text{if \ $d_1+d_2\geq n+1$}\\
\ \\
\left(\frac{d_2-1+n(d_1-1)}{(n-1)d_1}, \frac{d_1-1+n(d_2-1)}{(n-1)d_2}, 0\right),\ \left(1,1,\frac{n+1-(d_1+d_2)}{d_3}\right) & \text{if \ $d_1+d_2\leq n+1$}.\\
\end{cases}
\end{equation*}
\end{theorem}

\begin{theorem}\label{thm: main SSSS}{\rm{(Theorem \ref{thm: SSSS}, Example \ref{exa: SSSS})}}
Consider the log pair $(S_d, {S_{d_1}}|_{S_d}+{S_{d_2}}|_{S_d}+{S_{d_3}}|_{S_d})$, where $S_d, S_{d_1}, S_{d_2}, S_{d_3}$ are four general smooth hypersurfaces in $\bP^{n+1}$ of degrees $d, d_1, d_2, d_3$ with $n\geq 2$ and $d_1, d_2, d_3\leq n+2-d$. 
Suppose all the Fano complete intersections are K-semistable.
%Suppose  $S_d\cap S_{d_1}\cap S_{d_2}\cap S_{d_3}$ is K-semistable of dimension $n-3$ if $d_1+d_2+d_3\leq n+2-d$. 
\iffalse
Suppose the following two conditions:
\begin{enumerate}
\item $(S_d, {S_{d_1}}|_{S_d}+{S_{d_2}}|_{S_d}+{S_{d_3}}|_{S_d})$ is log canonical,
\item $S_d\cap S_{d_1}\cap S_{d_2}\cap S_{d_3}$ is K-semistable of dimension $n-3$ if $d_1+d_2+d_3\leq n+2-d$. 
\end{enumerate}
\fi
Then $\Kss(S_d, {S_{d_1}}|_{S_d}+{S_{d_2}}|_{S_d}+{S_{d_3}}|_{S_d})$ is a 
polytope $P$ generated by the following equations
\begin{equation*}
\begin{cases}
0\leq x\leq 1\\
0\leq y\leq 1\\
0\leq z\leq 1\\
d_1x+d_2y+d_3z\leq n+2-d\\
\frac{d_1(n+1)-(n+2-d)}{d_1(n+1)}-\frac{nx}{n+1}+\frac{d_2y}{d_1(n+1)}+\frac{d_3z}{d_1(n+1)}\geq 0\\
\frac{d_2(n+1)-(n+2-d)}{d_2(n+1)}-\frac{ny}{n+1}+\frac{d_1x}{d_2(n+1)}+\frac{d_3z}{d_2(n+1)}\geq 0\\
\frac{d_3(n+1)-(n+2-d)}{d_3(n+1)}-\frac{nz}{n+1}+\frac{d_1x}{d_3(n+1)}+\frac{d_2y}{d_3(n+1)}\geq 0.
\end{cases}
\end{equation*}
If $d_1+d_2+d_3<n+2-d$, then $P$ is generated by the extremal points 
$$(0,0,0), \  \left(\frac{(n+1)d_1-(n+2-d)}{nd_1}, 0,0\right),\  \left(0, \frac{(n+1)d_2-(n+2-d)}{nd_2}, 0\right), \  \left(0,0, \frac{(n+1)d_3-(n+2-d)}{nd_3}\right),$$

$$\left(0, \frac{d_3+nd_2-(n+2-d)}{(n-1)d_2}, \frac{d_2+nd_3-(n+2-d)}{(n-1)d_3}\right),$$  

$$\left(\frac{d_3+nd_1-(n+2-d)}{(n-1)d_1}, 0,  \frac{d_1+nd_3-(n+2-d)}{(n-1)d_3}\right), $$
 
$$\left(\frac{d_2+nd_1-(n+2-d)}{(n-1)d_1}, \frac{d_1+nd_2-(n+2-d)}{(n-1)d_2}, 0\right), $$

$$\left(\frac{d_2+d_3+(n-1)d_1-(n+2-d)}{(n-2)d_1}, \frac{d_1+d_3+(n-1)d_2-(n+2-d)}{(n-2)d_2}, \frac{d_1+d_2+(n-1)d_3-(n+2-d)}{(n-2)d_3}\right). $$
If $d_1+d_2+d_3\geq n+2-d$, then $P$ is generated by the extremal points
$$(0,0,0), \  \left(\frac{(n+1)d_1-(n+2-d)}{nd_1}, 0,0\right),\  \left(0, \frac{(n+1)d_2-(n+2-d)}{nd_2}, 0\right), \  \left(0,0, \frac{(n+1)d_3-(n+2-d)}{nd_3}\right),$$
and $A_1 ,B_1, A_2,B_2, A_3 ,B_3$, where
\begin{equation*}
A_1, B_1=
\begin{cases}
(0, 1, \frac{n+2-d-d_2}{d_3}),\ (0,\frac{n+2-d-d_3}{d_2}, 1) & \text{if \ $d_2+d_3\geq n+2-d$}\\
\ \\
\left(0, \frac{d_3+nd_2-(n+2-d)}{(n-1)d_2}, \frac{d_2+nd_3-(n+2-d)}{(n-1)d_3}\right),\ \left(\frac{n+2-d-(d_2+d_3)}{d_1},1,1\right) & \text{if \ $d_2+d_3\leq n+2-d$}\\
\end{cases}
\end{equation*}

\begin{equation*}
A_2, B_2=
\begin{cases}
(1, 0, \frac{n+2-d-d_1}{d_3}),\ (\frac{n+2-d-d_3}{d_1}, 0, 1) & \text{if \ $d_1+d_3\geq n+2-d$}\\
\ \\
\left(\frac{d_3+nd_1-(n+2-d)}{(n-1)d_1}, 0,  \frac{d_1+nd_3-(n+2-d)}{(n-1)d_3}\right),\ \left(1,\frac{n+2-d-(d_1+d_3)}{d_2},1\right) & \text{if \ $d_1+d_3\leq n+2-d$}\\
\end{cases}
\end{equation*}

\begin{equation*}
A_3, B_3=
\begin{cases}
(1, \frac{n+2-d-d_1}{d_2}, 0),\ (\frac{n+2-d-d_2}{d_1}, 1, 0) & \text{if \ $d_1+d_2\geq n+2-d$}\\
\ \\
\left(\frac{d_2+nd_1-(n+2-d)}{(n-1)d_1}, \frac{d_1+nd_2-(n+2-d)}{(n-1)d_2}, 0\right),\ \left(1,1, \frac{n+2-d-(d_1+d_2)}{d_3}\right) & \text{if \ $d_1+d_2\leq n+2-d$}.\\
\end{cases}
\end{equation*}
\end{theorem}

\begin{remark}
For all the results above, we pose the assumption that all the Fano complete intersections are K-semistable. Since we only treat general complete intersections, this condition is expected to be removed.
\end{remark}

We also present two even more general examples, which allow arbitrary numbers of components in the boundary. However, we do not plan to describe explicitly the extremal points of the K-semistable domains for them due to the complicated computation (but this could be done, even by hands). We do not list them here, please refer to Theorem \ref{thm: k components}, \ref{thm: k components'}.

\noindent
\subsection*{Acknowledgement}
The author is supported by the grant of European Research Council (ERC-804334).

\
\
\

\section{Preliminaries}\label{sec: preliminaries}

We say that $(X,\Delta)$ is a \emph{log pair} if $X$ is a normal projective variety and $\Delta$ is an effective $\bQ$-divisor on $X$ such that $K_X+\Delta$ is $\bQ$-Cartier.  The log pair $(X,\Delta)$ is called \emph{log Fano} if it admits klt singularities and $-(K_X+\Delta)$ is ample; if $\Delta=0$, we just say $X$ is a \emph{Fano variety}. 
The log pair $(X,\Delta)$ is called a \emph{log Calabi-Yau pair} if $K_X+\Delta\sim_\bQ 0$. 
For various types of singularities in birational geometry, e.g.  klt and lc singularities, we refer to \cite{KM98,Kollar13}.

\subsection{K-stability}

Let $(X,\Delta)$ be a log pair. Suppose $f\colon Y\to X$ is a proper birational morphism between normal varieties and $E$ is a prime divisor on $Y$, we say that $E$ is a prime divisor over $X$ and define the following invariant
$$A_{X,\Delta}(E):=1+\ord_E(K_Y-f^*(K_X+\Delta)), $$
which is called the \emph{log discrepancy} of $E$ associated to the log pair $(X,\Delta)$.
If $(X,\Delta)$ is a log Fano pair, we define the following invariant
$$S_{X,\Delta}(E):=\frac{1}{\vol(-K_X-\Delta)}\int_0^\infty \vol(-f^*(K_X+\Delta)-tE){\rm{d}}t .$$
Put $\beta_{X,\Delta}(E):=A_{X,\Delta}(E)-S_{X,\Delta}(E)$. By the works \cite{Fuj19, Li17}, one can define K-stability of a log Fano pair by beta criterion as follows.
\begin{definition}\label{def: kss}
Let $(X,\Delta)$ be a log Fano pair. 
We say that $(X,\Delta)$ is \emph{K-semistable} if $\beta_{X,\Delta}(E)\geq 0$ for any prime divisor $E$ over $X$.
\end{definition}

By the works \cite{Oda13, BHJ17}, one can define K-stability of a log Calabi-Yau pair by posing singularity condition.

\begin{definition}
Let $(X, \Delta)$ be a log Calabi-Yau pair, i.e. $K_X+\Delta\sim_\bQ 0$. We say  $(X,\Delta)$ is \emph{K-semistable} if  $(X, \Delta)$ is log canonical.
\end{definition}

\begin{remark}\label{rem: kss hypersurface}
It is a hard problem to test K-stability for explicit Fano varieties, and people have made some progress on several special examples in the past few years, e.g. \cite{AZ22, AZ23}. However, in general setting, with the help of theoretical development on algebraic K-stability (\cite{BLX22, Xu20}), we at least know
\begin{enumerate}
\item A general hypersurface $S_d\subset \bP^{n+1}$ of degree $d$ with $d\leq n+2$ is K-semistable.
\item The complete intersection of $m$ general hypersurfaces $S_{i}\subset \bP^{n+m}$ of degree $d$ ($ 1\leq i\leq m$) with $md\leq n+m+1$ is K-semistable (e.g. \cite{AGP06}).
\end{enumerate}
\end{remark}

\subsection{K-semistable domain}

In \cite{LZ23}, a set $\mE$ of log pairs is defined. Fix two positive integers $d$ and $k$, a positive number $v$,  and a finite set $I$ of non-negative rational numbers, we consider the set $\mE:=\mE(d,k,v, I)$ of log pairs $(X, \sum_{i=1}^k D_i)$ satisfying the following conditions:
\begin{enumerate}
\item $X$ is a Fano variety of dimension $d$ and $(-K_X)^d=v$;
\item $0\leq D_i\sim_\bQ -K_X$ for every $1\leq i\leq k$;
\item the coefficients of $D_i$ are contained in $I$;
\item there exists $(c_1,...,c_k)\in \Delta^k$ such that $(X, \sum_i c_iD_i)$ is K-semistable, where $\Delta^k:=\{(c_1,...,c_k)\ |\ \textit{$c_i\in [0,1)\cap \bQ$ and $0\leq \sum_i c_i<1$}\}$.
\end{enumerate}

\begin{definition}\label{def: kss domain}
Let $(X, \sum_{i=1}^kD_i)$ be a log pair in the set $\mE$. We define the \emph{K-semistable domain} of $(X, \sum_{i=1}^kD_i)$ as follows:
$$\mathrm{Kss}(X, \sum_{i=1}^kD_i):=\overline{\{(c_1,...,c_k)\ |\ \textit{$c_i\in [0,1)\cap \bQ$ and $(X, \sum_{i=1}^kc_iD_i)$ is K-semistable}\}}. $$
The overline in the definition means taking the closure.
\end{definition}

\begin{remark}\label{rem: kss domain}
In the definition of the set $\mE$, we just assume $D_i\sim_\bQ -K_X$ for convenience, since under this assumption the K-semistable domain $\Kss(X, \sum_{i=1}^kD_i)$ lies in the simplex $\overline{\Delta^k}$. If we replace $D_i\sim_\bQ -K_X$ with that $D_i$ is proportional to $-K_X$, then clearly one could define the K-semistable domain by the same way. The proportional condition allows us to apply the following interpolation property for K-stability, which we will use frequently: if $(X, \Delta_1)$ and $(X, \Delta_2)$ are both K-semistable log pairs (log Fano or log Calabi-Yau), where $\Delta_i$ are propotional to $-K_X$, then $(X, t\Delta_1+(1-t)\Delta_2)$ is also K-semistable for any $t\in [0,1]\cap \bQ$.
\end{remark}

Let $V$ be a Fano manifold of dimension $n$, and $S$ a smooth divisor on $V$ such that $S\sim_\bQ -\lambda K_V$ for some positive rational number $\lambda$. Recall that
$$\mathrm{Kss}(V,S)=\overline{\{a\in [0,1)\cap \bQ\ |\ \textit{$(V,aS)$ is K-semistable} \}}. $$

\begin{lemma}{\rm{(\cite{ZZ22})}}\label{lem: zz general}
Notation as above, suppose $V$ and S are both K-semistable and $0<\lambda\leq 1$, then $\mathrm{Kss}(V,S)= [0,1-\frac{r}{n}]$, where $r=\frac{1}{\lambda}-1$. 
\end{lemma}

As a special case, we have

\begin{lemma}\label{lem: zz}
Let $(\bP^n, S_d)$ be a log pair where $S_d$ is a general smooth hypersurface of degree $d\leq n+1$.  Then we have $\mathrm{Kss}(\bP^n, S_d)=[0, 1-\frac{r}{n}]$, where $r=\frac{n+1-d}{d}$.
\end{lemma}

The following lemma is also well known, see e.g. \cite[Prop 2.11]{LZ22}.
\begin{lemma}\label{lem:cone stability}
Let $(V,\Delta)$ be an n-dimensional log Fano pair, and L an ample line bundle on V such that $L\sim_\bQ -\frac{1}{r}(K_V+\Delta)$ for some $0<r\leq n+1$. Suppose Y is the projective cone over V associated to L with infinite divisor $V_\infty$, then $(V,\Delta)$ is K-semistable (resp. K-polystable) if and only if $(Y,\Delta_Y+(1-\frac{r}{n+1})V_\infty)$ is K-semistable (resp. K-polystable), where $\Delta_Y$ is the divisor on Y naturally extended by $\Delta$.
\end{lemma}

\begin{remark}\label{rem: facts}
Applying Lemma \ref{lem: zz general}, a simple computation tells us the following facts:

\begin{enumerate}
\item Let $S_d\subset \bP^n$ be a general smooth hypersurface of degree $d\leq n+1$. Then $(\bP^n, tS_d)$ is K-semistable for any $t\in [0, \frac{(n+1)(d-1)}{nd}]$. 
\item Let $S_d, S_{d'}\subset \bP^{n+1}$ be two general smooth hypersurfaces of degrees $d, d'$ with $d+d'\leq n+2$. Suppose $S_d\cap S_{d'}$ is K-semistable, then $(S_d, t{S_{d'}}|_{S_d})$ is K-semistable  for any $t\in [0, 1-\frac{n+2-d-d'}{nd'}]$.
\item Let $S_d, S_{d_1}, S_{d_2}\subset \bP^{n+2}$ be three general smooth hypersurfaces of degrees $d, d_1, d_2$ with $d_1+d_2\leq n+3-d$. Suppose $S_d\cap S_{d_1}$ and $S_d\cap S_{d_1}\cap S_{d_2}$ are K-semistable,
then 
$$\left({S_{d_1}}|_{S_d}, t({S_{d_2}}|_{S_d})|_{{S_{d_1}}|_{S_d}}\right)$$
is K-semistable  for any $t\in [0, 1-\frac{n+3-d-d_1-d_2}{nd_2}]$.
\item Let $S_d, S_{d_1}, S_{d_2}, S_{d_3}\subset \bP^{n+3}$ be four general smooth hypersurfaces of degrees $d, d_1, d_2, d_3$ with $d_1+d_2+d_3\leq n+4-d$.
Suppose $S_d\cap S_{d_1}\cap S_{d_2}$ and $S_d\cap S_{d_1}\cap S_{d_2}\cap S_{d_3}$ are K-semistable,
then 
$$\left( ({S_{d_2}}|_{S_d})|_{{S_{d_1}}|_{S_d}}, t\left(({S_{d_3}}|_{S_d})|_{{S_{d_1}}|_{S_d}}\right)|_{({S_{d_2}}|_{S_d})|_{{S_{d_1}}|_{S_d}}}\right)=(S_{d_2}\cap S_{d_1}\cap S_d, t S_{d_3}\cap S_{d_2}\cap S_{d_1}\cap S_d )$$ 
is K-semistable  for any $t\in [0, 1-\frac{n+4-d-d_1-d_2-d_3}{nd_3}]$.

\end{enumerate}
\end{remark}

\iffalse
\subsection{K-semistable degeneration}

Let $(X, \Delta)$ be a log Fano pair, then the delta invariant of $(X, \Delta)$ is defined as 
$$\delta(X, \Delta):=\inf_{E} \frac{A_{X, \Delta}(E)}{S_{X, \Delta}(E)}, $$
where $E$ runs through all prime divisors over $X$ (e.g. \cite{FO18, BJ20}). By Definition \ref{def: kss} we see that $(X, \Delta)$ is K-semistable if and only if $\delta(X, \Delta)\geq 1$.

\begin{theorem}{\rm{(\cite{LXZ22})}}\label{thm: lxz22}
Let $(X, \Delta)$ be a log Fano pair with $\delta(X, \Delta)\leq 1$, then there exists a prime divisor $E$ over $X$ such that $\delta(X, \Delta)=\frac{A_{X, \Delta}(E)}{S_{X, \Delta}(E)}$.
\end{theorem}

\begin{theorem}{\rm{(\cite{BLZ22})}}\label{thm: blz22}
Let $(X, \Delta)$ be a log Fano pair with $\delta(X, \Delta)\leq 1$, and $E$ a prime divisor over $X$ such that $\delta(X, \Delta)=\frac{A_{X, \Delta}(E)}{S_{X, \Delta}(E)}$. Then $E$ induces a special test configuration $(\mX, \Delta_\tc)\to \bA^1$ of $(X, \Delta)$ such that $\delta(\mX_0, \Delta_{\tc, 0})=\delta(X, \Delta)$. In particular, if $(X, \Delta)$ is K-semistable, then the central fiber $(\mX_0, \Delta_{\tc, 0})$ is also K-semistable.
\end{theorem}
\fi

\begin{remark}\label{rem: kss degeneration}
In this paper, we are interested in the K-semistable degeneration of a log Fano pair $(V, \Delta)$. Let $S$ be a prime Cartier divisor on $V$ with $S\sim_\bQ -\lambda (K_V+\Delta)$ for some rational $0<\lambda<1$. By the finite generation of the following graded ring:
$$\bigoplus_{m\in \bN}\bigoplus_{j\in \bN}H^0(V, -mr(K_V+\Delta)-jS), $$
where $r$ is a fixed divisible positive integer, one could construct a test configuration of $(V, \Delta)$, denoted by $(\mV, \Delta_{\tc})\to \bA^1$, such that the central fiber $\mV_0$ is isomorphic to the projective cone over $S$ with respect to the polarization $S|_S\sim_\bQ -\frac{1}{r}K_S$ for $r=\frac{1}{\lambda}-1$, and the infinite divisor of $\mV_0$ is isomorphic to $S$. \end{remark}

\
\
\

\section{$(\bP^n, S_{d_1}+S_{d_2})$}

In this section, we compute $\Kss(\bP^n, S_{d_1}+S_{d_2})$, where $S_{d_1}, S_{d_2}$ are two general smooth hypersurfaces in $\bP^n$ of degrees $d_1, d_2$ with $d_1, d_2\leq n+1$ and $n\geq 2$.

\begin{theorem}\label{thm: PSS}
Let $n, d_1, d_2$ be three positive integers satisfying $n\geq 2$ and $d_1+d_2<n+1$. Suppose $S_{d_1}\cap S_{d_2}$ is K-semistable.
Then the log Fano pair 
$$\left(\bP^n, \frac{d_2-1+n(d_1-1)}{(n-1)d_1}S_{d_1}+ \frac{d_1-1+n(d_2-1)}{(n-1)d_2}S_{d_2}\right) $$
is K-semistable.
\end{theorem}

\begin{proof}

By Remark \ref{rem: kss degeneration},  $S_{d_1}$ induces a test configuration of the log Fano pair 
$$\left(\bP^n, \frac{d_2-1+n(d_1-1)}{(n-1)d_1}S_{d_1}+ \frac{d_1-1+n(d_2-1)}{(n-1)d_2}S_{d_2}\right) $$
such that the central fiber, denoted by
$$\left(X, \frac{d_2-1+n(d_1-1)}{(n-1)d_1}D+ \frac{d_1-1+n(d_2-1)}{(n-1)d_2}D'\right) $$
is the projective cone over 
$$\left(S_{d_1}, \frac{d_1-1+n(d_2-1)}{(n-1)d_2}{S_{d_2}}|_{S_{d_1}}\right)$$ 
with respect to the polarization $\mO_{S_{d_1}}(d_1):=i^*\mO_{\bP^{n+1}}(d_1)$,  where $D\cong S_{d_1}$ is the infinite divisor and $i: {S_{d_1}}\to \bP^{n+1}$ is the natural embedding.
It suffices to show that the log pair 
$$\left(X, \frac{d_2-1+n(d_1-1)}{(n-1)d_1}D+ \frac{d_1-1+n(d_2-1)}{(n-1)d_2}D'\right) $$
is a K-semistable log Fano pair.  
We have the following computation:
$$\mO_{S_{d_1}}(d_1)\sim_\bQ -\frac{1}{r}\left(K_{S_{d_1}}+\frac{d_1-1+n(d_2-1)}{(n-1)d_2}{S_{d_2}}|_{S_{d_1}}\right)$$
and  
$$\frac{d_2-1+n(d_1-1)}{(n-1)d_1}=1-\frac{r}{n}$$
for 
$$r= \frac{n+1-d_1-\frac{d_1-1+n(d_2-1)}{n-1}}{d_1}.$$ 
Note that $(S_{d_1}, a{S_{d_2}}|_{S_{d_1}})$ is K-semistable for $a\in [0, 1-\frac{n+1-d_1-d_2}{(n-1)d_2}]$ by Remark \ref{rem: facts}, thus
$$\left(S_{d_1}, \frac{d_1-1+n(d_2-1)}{(n-1)d_2}{S_{d_2}}|_{S_{d_1}}\right)$$
is K-semistable. Applying  Lemma \ref{lem:cone stability}, we see that $$\left(X, \frac{d_2-1+n(d_1-1)}{(n-1)d_1}D+ \frac{d_1-1+n(d_2-1)}{(n-1)d_2}D'\right) $$
 is also K-semistable.
\end{proof}

\begin{example}\label{exa: PSS}
Consider the log pair $(\bP^n, S_{d_1}+S_{d_2})$, where $S_{d_1}, S_{d_2}$ are two general smooth hypersurfaces of degrees $d_1, d_2\leq n+1$ and $n\geq 2$. 
\iffalse 
Suppose the following two conditions:
\begin{enumerate}
\item $(\bP^n, S_{d_1}+S_{d_2})$ is log canonical,
\item $S_{d_1}\cap S_{d_2}$ is K-semistable of dimension $n-2$ if $d_1+d_2\leq n+1$. 
\end{enumerate}
\fi
Suppose $S_{d_1}\cap S_{d_2}$ is K-semistable if $d_1+d_2\leq n+1$. We want to compute $\Kss(\bP^n, S_{d_1}+S_{d_2})$.

Suppose $(x,y)\in \Kss(\bP^n, S_{d_1}+S_{d_2})$, then $(\bP^n, xS_{d_1}+yS_{d_2})$ is K-semistable. Applying beta criterion, we have
\begin{align*}
\beta_{\bP^n, xS_{d_1}+yS_{d_2}}(S_{d_1})\ &=\ 1-x-\frac{1}{(n+1-d_1x-d_2y)^n}\int_0^\frac{{n+1-d_1x-d_2y}}{d_1}(n+1-d_1x-d_2y-d_1t)^n\dt \\
&=\ 1-x-\frac{n+1-d_1x-d_2y}{d_1(n+1)}\\
&=\ \frac{d_1-1}{d_1}-\frac{nx}{n+1}+\frac{d_2y}{d_1(n+1)}\geq 0,
\end{align*}
\begin{align*}
\beta_{\bP^n, xS_{d_1}+yS_{d_2}}(S_{d_2})\ &=\ 1-y-\frac{1}{(n+1-d_1x-d_2y)^n}\int_0^\frac{{n+1-d_1x-d_2y}}{d_2}(n+1-d_1x-d_2y-d_2t)^n\dt \\
&=\ 1-y-\frac{n+1-d_1x-d_2y}{d_2(n+1)}\\
&=\ \frac{d_2-1}{d_2}-\frac{ny}{n+1}+\frac{d_1x}{d_2(n+1)}\geq 0.
\end{align*}
We denote by $P$ the polytope generated by the following equations
\begin{equation*}
\begin{cases}
0\leq x\leq 1\\
0\leq y\leq 1\\
d_1x+d_2y\leq n+1\\
\frac{d_1-1}{d_1}-\frac{nx}{n+1}+\frac{d_2y}{d_1(n+1)}\geq 0\\
\frac{d_2-1}{d_2}-\frac{ny}{n+1}+\frac{d_1x}{d_2(n+1)}\geq 0.\\
\end{cases}
\end{equation*}
If $d_1+d_2<n+1$, then $P$ is generated by the extremal points 
$$(0,0), \ \ \left(\frac{(n+1)(d_1-1)}{nd_1}, 0\right),\ \ \left(0, \frac{(n+1)(d_2-1)}{nd_2}\right), \ \ \left(\frac{d_2-1+n(d_1-1)}{(n-1)d_1}, \frac{d_1-1+n(d_2-1)}{(n-1)d_2}\right). $$
If $d_1+d_2\geq n+1$, then $P$ is generated by the extremal points
$$(0,0), \ \ \left(\frac{(n+1)(d_1-1)}{nd_1}, 0\right),\ \ \left(0, \frac{(n+1)(d_2-1)}{nd_2}\right), \ \ \left(1, \frac{n+1-d_1}{d_2}\right), \ \ \left(\frac{n+1-d_2}{d_1}, 1\right). $$ 
In either case, these points correspond to K-semistable log pairs by Theorem \ref{thm: PSS} and Remark \ref{rem: facts}, thus we exactly have 
$\Kss(\bP^n, S_{d_1}+S_{d_2})=P$.
\end{example}

The above example in fact generalizes \cite[Theorem 1.4, 1.5]{LZ23}.

\begin{example}
In Example \ref{exa: PSS}, take $n\geq 2, d_1=2, d_2=1$. We denote $S_{d_1}, S_{d_2}$ by $Q, L$. Then $\Kss(\bP^n, Q+L)$ is generated by the extremal points
$$(0, 0),\ \ \left(\frac{n+1}{2n}, 0\right),\ \ \left(\frac{n}{2(n-1)}, \frac{1}{n-1}\right). $$
This is exactly \cite[Theorem 1.4]{LZ23}.
\end{example}

\begin{example}
In Example \ref{exa: PSS}, take $n\geq 3, d_1=2, d_2=2$. We denote $S_{d_1}, S_{d_2}$ by $Q, Q'$. Then $\Kss(\bP^n, Q+Q')$ is generated by the extremal points
$$(0, 0),\ \ \left(\frac{n+1}{2n}, 0\right),\ \ \left(0, \frac{n+1}{2n}\right),\ \ \left(\frac{n+1}{2(n-1)}, \frac{n+1}{2(n-1)}\right). $$
This is exactly \cite[Theorem 1.5]{LZ23}.
\end{example}

\
\
\

\section{$(S_d, {S_{d_1}}|_{S_d}+{S_{d_2}}|_{S_d})$}

In this section, we compute $\Kss(S_d, {S_{d_1}}|_{S_d}+{S_{d_2}}|_{S_d})$, where $S_d, S_{d_1}, S_{d_2}$ are three general smooth hypersurfaces in $\bP^{n+1}$ of degrees $d, d_1, d_2$ with $d_1, d_2\leq n+2-d$ and $n\geq 2$.

\begin{theorem}\label{thm: SSS}
Let $d, d_1, d_2$ be three positive integers satisfying $d_1+d_2<n+2-d$ and $n\geq 2$. 
Suppose all the Fano complete intersections are K-semistable.
%Suppose $S_d\cap S_{d_1}\cap S_{d_2}$ is K-semistable of dimension $n-2$.
Then the log Fano pair
$$\left(S_d, \frac{d_2+nd_1-(n+2-d)}{(n-1)d_1}{S_{d_1}}|_{S_d} +\frac{d_1+nd_2-(n+2-d)}{(n-1)d_2}{S_{d_2}}|_{S_d}\right)$$
is K-semistable.
\end{theorem}

\begin{proof}
By Remark \ref{rem: kss degeneration},  $S_d\cap S_{d_1}$ induces a test configuration of the log Fano pair 
$$\left(S_d, \frac{d_2+nd_1-(n+2-d)}{(n-1)d_1}{S_{d_1}}|_{S_d} +\frac{d_1+nd_2-(n+2-d)}{(n-1)d_2}{S_{d_2}}|_{S_d}\right)$$
such that the central fiber, denoted by
$$\left(X, \frac{d_2+nd_1-(n+2-d)}{(n-1)d_1}D +\frac{d_1+nd_2-(n+2-d)}{(n-1)d_2}D'\right),$$ 
is the projective cone over 
$$(S_d\cap S_{d_1}, \frac{d_1+nd_2-(n+2-d)}{(n-1)d_2}S_d\cap S_{d_1}\cap S_{d_2})=\left({S_{d_1}}|_{S_d}, \frac{d_1+nd_2-(n+2-d)}{(n-1)d_2}({S_{d_2}}|_{S_d})|_{{S_{d_1}}|_{S_d}}\right)$$ 
with respect to the polarization $\mO_{{S_{d_1}}|_{S_d}}(d_1):=i^*\mO_{\bP^{n+1}}(d_1)$,  where $D\cong S_d\cap S_{d_1}$ is the infinite divisor and $i: {S_{d_1}}|_{S_d}\to \bP^{n+1}$ is the natural embedding.
It suffices to show that the log pair 
$$\left(X, \frac{d_2+nd_1-(n+2-d)}{(n-1)d_1}D +\frac{d_1+nd_2-(n+2-d)}{(n-1)d_2}D'\right)$$ 
is a K-semistable log Fano pair.  We have the following computation:
$$\mO_{{S_{d_1}}|_{S_d}}(d_1)\sim_\bQ -\frac{1}{r}\left(K_{{S_{d_1}}|_{S_d}}+\frac{d_1+nd_2-(n+2-d)}{(n-1)d_2}({S_{d_2}}|_{S_d})|_{{S_{d_1}}|_{S_d}}\right)$$
and  
$$\frac{d_2+nd_1-(n+2-d)}{(n-1)d_1}=1-\frac{r}{n}$$
for 
$$r= \frac{n+2-d-d_1-\frac{d_1+nd_2-(n+2-d)}{n-1}}{d_1}.$$ 
By Lemma \ref{lem:cone stability}, it suffices to show that 
$$(S_d\cap S_{d_1}, \frac{d_1+nd_2-(n+2-d)}{(n-1)d_2}S_d\cap S_{d_1}\cap S_{d_2})$$ 
is K-semistable, which is implied by Lemma \ref{lem: zz general} or Remark \ref{rem: facts}.
\iffalse
By Remark \ref{rem: kss degeneration}, $S_d\cap S_{d_1}\cap S_{d_2}$ induces a test configuration of the above log Fano pair such that the central fiber, denoted by $(Y, \frac{d_1+nd_2-(n+2-d)}{(n-1)d_2}B)$, is the projective cone over $S_d\cap S_{d_1}\cap S_{d_2}$ with respect to the polarization $\mO_{S_d\cap S_{d_1}\cap S_{d_2}}(d_2)$, and $B\cong S_d\cap S_{d_1}\cap S_{d_2}$ is the infinite divisor. It suffices to show that $(Y, \frac{d_1+nd_2-(n+2-d)}{(n-1)d_2}B)$ is K-semistable.
We have the following computation 
$$\mO_{S_d\cap S_{d_1}\cap S_{d_2}}(d_2)\sim_\bQ -\frac{1}{r}K_{S_d\cap S_{d_1}\cap S_{d_2}}$$
and  
$$\frac{d_1+nd_2-(n+2-d)}{(n-1)d_2}=1-\frac{r}{n-1}$$
for 
$$r= \frac{n+2-d-d_1-d_2}{d_2}.$$ 
Since $S_d\cap S_{d_1}\cap S_{d_2}$ is K-semistable, applying  Lemma \ref{lem:cone stability} again we see that the log Fano pair $(Y, \frac{d_1+nd_2-(n+2-d)}{(n-1)d_2}B)$ is K-semistable.
\fi
\end{proof}

\begin{example}\label{exa: SSS}
Consider the log pair $(S_d, {S_{d_1}}|_{S_d}+{S_{d_2}}|_{S_d})$, where $S_d, S_{d_1}, S_{d_2}$ are three general smooth hypersurfaces in $\bP^{n+1}$ of degrees $d, d_1, d_2$ with $n\geq 2$ and $ d_1, d_2\leq n+2-d$. 
\iffalse
Suppose the following two conditions:
\begin{enumerate}
\item $(S_d, {S_{d_1}}|_{S_d}+{S_{d_2}}|_{S_d})$ is log canonical,
\item $S_d\cap S_{d_1}\cap S_{d_2}$ is K-semistable of dimension $n-2$ if $d_1+d_2\leq n+2-d$. 
\end{enumerate}
\fi
Suppose all the Fano complete intersections are K-semistable.
%Suppose $S_d\cap S_{d_1}\cap S_{d_2}$ is K-semistable of dimension $n-2$ if $d_1+d_2\leq n+2-d$. 
We want to compute $\Kss(S_d, {S_{d_1}}|_{S_d}+{S_{d_2}}|_{S_d})$.

Suppose $(x,y)\in \Kss(S_d, {S_{d_1}}|_{S_d}+{S_{d_2}}|_{S_d})$, then $\Kss(S_d, x{S_{d_1}}|_{S_d}+y{S_{d_2}}|_{S_d})$ is K-semistable. Applying beta criterion, we have
\begin{align*}
&\ \beta_{S_d, x{S_{d_1}}|_{S_d}+y{S_{d_2}}|_{S_d}}({S_{d_1}}|_{S_d})\\
=&\ 1-x-\frac{1}{d(n+2-d-d_1x-d_2y)^n}\int_0^\frac{{n+2-d-d_1x-d_2y}}{d_1}d(n+2-d-d_1x-d_2y-d_1t)^n\dt \\
=&\ 1-x-\frac{n+2-d-d_1x-d_2y}{d_1(n+1)}\\
=&\ \frac{d_1(n+1)-(n+2-d)}{d_1(n+1)}-\frac{nx}{n+1}+\frac{d_2y}{d_1(n+1)}\geq 0,
\end{align*}
\begin{align*}
&\ \beta_{S_d, x{S_{d_1}}|_{S_d}+y{S_{d_2}}|_{S_d}}({S_{d_2}}|_{S_d})\\
=&\ 1-y-\frac{1}{d(n+2-d-d_1x-d_2y)^n}\int_0^\frac{{n+2-d-d_1x-d_2y}}{d_2}d(n+2-d-d_1x-d_2y-d_2t)^n\dt \\
=&\ 1-y-\frac{n+2-d-d_1x-d_2y}{d_2(n+1)}\\
=&\ \frac{d_2(n+1)-(n+2-d)}{d_2(n+1)}-\frac{ny}{n+1}+\frac{d_1x}{d_2(n+1)}\geq 0.
\end{align*}
We denote by $P$ the polytope generated by the following equations
\begin{equation*}
\begin{cases}
0\leq x\leq 1\\
0\leq y\leq 1\\
d_1x+d_2y\leq n+2-d\\
\frac{d_1(n+1)-(n+2-d)}{d_1(n+1)}-\frac{nx}{n+1}+\frac{d_2y}{d_1(n+1)}\geq 0\\
\frac{d_2(n+1)-(n+2-d)}{d_2(n+1)}-\frac{ny}{n+1}+\frac{d_1x}{d_2(n+1)}\geq 0.\\
\end{cases}
\end{equation*}
If $d_1+d_2<n+2-d$, then $P$ is generated by the extremal points 
$$(0,0), \ \ \left(\frac{(n+1)d_1-(n+2-d)}{nd_1}, 0\right),\ \ \left(0, \frac{(n+1)d_2-(n+2-d)}{nd_2}\right)$$
and
$$\left(\frac{d_2+nd_1-(n+2-d)}{(n-1)d_1}, \frac{d_1+nd_2-(n+2-d)}{(n-1)d_2}\right). $$
If $d_1+d_2\geq n+2-d$, then $P$ is generated by the extremal points
$$(0,0), \ \ \left(\frac{(n+1)d_1-(n+2-d)}{nd_1}, 0\right),\ \ \left(0, \frac{(n+1)d_2-(n+2-d)}{nd_2}\right)$$
and
$$\left(1, \frac{n+2-d-d_1}{d_2}\right), \ \ \left(\frac{n+2-d-d_2}{d_1}, 1\right). $$ 
In either case, these points correspond to K-semistable log pairs by Theorem \ref{thm: SSS} and Remark \ref{rem: facts}, thus we exactly have 
$\Kss(S_d, {S_{d_1}}|_{S_d}+{S_{d_2}}|_{S_d})=P$.
\end{example}

\begin{example}
In Example \ref{exa: SSS}, take $n\geq 2, d=2, d_1=1, d_2=1$. Denote $S_d, S_{d_1}, S_{d_2}$ by $Q, H_1, H_2$. Then $\Kss(Q, {H_1}|_Q+{H_2}|_Q)$ is generated by the extremal points
$$(0,0),\ \ (\frac{1}{n}, 0),\ \ (0, \frac{1}{n}),\ \ (\frac{1}{n-1}, \frac{1}{n-1}).  $$
When $n=3$, this example is computed by \cite[Theorem 1.1 (2)]{Loginov23}.

\end{example}

\
\
\

\section{$({S_{d_1}}|_{S_d}, {({S_{d_2}}|_{S_d})}|_{{S_{d_1}}|_{S_d}}+{({S_{d_3}}|_{S_d})}|_{{S_{d_1}}|_{S_d}})$}

In this section, we want to compute
$\Kss\left({S_{d_1}}|_{S_d}, {({S_{d_2}}|_{S_d})}|_{{S_{d_1}}|_{S_d}}+{({S_{d_3}}|_{S_d})}|_{{S_{d_1}}|_{S_d}}\right),$ where $S_d, S_{d_1}, S_{d_2}, S_{d_3}$ are general smooth hypersurfaces in $\bP^{n+2}$ of degrees $d, d_1, d_2, d_3$ with $d_2, d_3\leq n+3-d-d_1$ and $n\geq 2$.

\begin{theorem}\label{thm: SSS|S}
Let $d, d_1, d_2, d_3$ be four positive integers satisfying $d_1+d_2+d_3<n+3-d$ and $n\geq 2$. 
Suppose all the Fano complete intersections are K-semistable.
%Suppose $S_d\cap S_{d_1}\cap S_{d_2}\cap S_{d_3}$ is K-semistable of dimension $n-2$. 
Then the log Fano pair
$$\left({S_{d_1}}|_{S_d}, b{({S_{d_2}}|_{S_d})}|_{{S_{d_1}}|_{S_d}}+c{({S_{d_3}}|_{S_d})}|_{{S_{d_1}}|_{S_d}}\right)$$
is K-semistable for 
$$b=\frac{d_1+d_3+nd_2-(n+3-d)}{(n-1)d_2},$$
$$c=\frac{d_1+d_2+nd_3-(n+3-d)}{(n-1)d_3}. $$
\end{theorem}

\begin{proof}
By Remark \ref{rem: kss degeneration},  $S_d\cap S_{d_1}\cap S_{d_2}$ induces a test configuration of 
$$(S_d\cap S_{d_1}, bS_d\cap S_{d_1}\cap S_{d_2}+cS_d\cap S_{d_1}\cap S_{d_3})=\left({S_{d_1}}|_{S_d}, b{({S_{d_2}}|_{S_d})}|_{{S_{d_1}}|_{S_d}}+c{({S_{d_3}}|_{S_d})}|_{{S_{d_1}}|_{S_d}}\right)$$
such that the central fiber, denoted by $(X, bD+cD')$, is the projective cone over 
$$\left( ({S_{d_2}}|_{S_d})|_{{S_{d_1}}|_{S_d}}, c\left(({S_{d_3}}|_{S_d})|_{{S_{d_1}}|_{S_d}}\right)|_{\left(({S_{d_2}}|_{S_d})|_{{S_{d_1}}|_{S_d}}\right)}\right)=(S_{d_2}\cap S_{d_1}\cap S_d, c S_{d_3}\cap S_{d_2}\cap S_{d_1}\cap S_d )$$ 
with respect to the polarization $\mO_{S_{d_2}\cap S_{d_1}\cap S_d}(d_2)$ and $D\cong S_d\cap S_{d_1}\cap S_{d_2}$ is the infinite divisor. It suffices to show that $(X, bD+cD')$ is a K-semistable log Fano pair.
We have the following computation:
$$\mO_{S_{d_2}\cap S_{d_1}\cap S_d}(d_2)\sim_\bQ -\frac{1}{r}(K_{S_{d_2}\cap S_{d_1}\cap S_d}+cS_{d_3}\cap S_{d_2}\cap S_{d_1}\cap S_d)$$
and  
$$b=1-\frac{r}{n}$$
for 
$$r= \frac{n+3-d-d_1-d_2-cd_3}{d_2}.$$ 
Note that 
$$\left( ({S_{d_2}}|_{S_d})|_{{S_{d_1}}|_{S_d}}, t\left(({S_{d_3}}|_{S_d})|_{{S_{d_1}}|_{S_d}}\right)|_{\left(({S_{d_2}}|_{S_d})|_{{S_{d_1}}|_{S_d}}\right)}\right)$$ 
is K-semistable for $t\in [0, 1-\frac{n+3-d-d_1-d_2-d_3}{(n-1)d_3}]$ by Remark \ref{rem: facts}, thus
$$\left( ({S_{d_2}}|_{S_d})|_{{S_{d_1}}|_{S_d}}, c\left(({S_{d_3}}|_{S_d})|_{{S_{d_1}}|_{S_d}}\right)|_{\left(({S_{d_2}}|_{S_d})|_{{S_{d_1}}|_{S_d}}\right)}\right)$$ 
is K-semistable. Applying  Lemma \ref{lem:cone stability}, we see that 
$(X, bD +cD')$
 is also K-semistable.
\end{proof}

\begin{example}\label{exa: SSS|S}
Consider the log pair $\left({S_{d_1}}|_{S_d}, {({S_{d_2}}|_{S_d})}|_{{S_{d_1}}|_{S_d}}+{({S_{d_3}}|_{S_d})}|_{{S_{d_1}}|_{S_d}}\right)$, where $S_d, S_{d_1}, S_{d_2}$ and $S_{d_3}$ are four general smooth hypersurfaces of degrees $ d, d_1, d_2, d_3$ in $\bP^{n+2}$ with $n\geq 2$ and $d_2, d_3\leq n+3-d-d_1$.
\iffalse
Suppose the following two conditions:
\begin{enumerate}
\item $\left({S_{d_1}}|_{S_d}, {({S_{d_2}}|_{S_d})}|_{{S_{d_1}}|_{S_d}}+{({S_{d_3}}|_{S_d})}|_{{S_{d_1}}|_{S_d}}\right)$ is log canonical,
\item $S_d\cap S_{d_1}\cap S_{d_2}\cap S_{d_3}$ is K-semistable of dimension $n-2$ if $d_2+d_3\leq n+3-d-d_1$. 
\end{enumerate}
\fi
Suppose all the Fano complete intersections are K-semistable.
%Suppose $S_d\cap S_{d_1}\cap S_{d_2}\cap S_{d_3}$ is K-semistable of dimension $n-2$ if $d_2+d_3\leq n+3-d-d_1$. 
We want to compute 
$$\Kss\left({S_{d_1}}|_{S_d}, {({S_{d_2}}|_{S_d})}|_{{S_{d_1}}|_{S_d}}+{({S_{d_3}}|_{S_d})}|_{{S_{d_1}}|_{S_d}}\right).$$

Suppose $(x,y)\in \Kss\left({S_{d_1}}|_{S_d}, {({S_{d_2}}|_{S_d})}|_{{S_{d_1}}|_{S_d}}+{({S_{d_3}}|_{S_d})}|_{{S_{d_1}}|_{S_d}}\right)$, then 
$$\left({S_{d_1}}|_{S_d}, x{({S_{d_2}}|_{S_d})}|_{{S_{d_1}}|_{S_d}}+y{({S_{d_3}}|_{S_d})}|_{{S_{d_1}}|_{S_d}}\right)$$ is K-semistable. Applying beta criterion, we have
\begin{align*}
&\ \beta_{{S_{d_1}}|_{S_d}, x{({S_{d_2}}|_{S_d})}|_{{S_{d_1}}|_{S_d}}+y{({S_{d_3}}|_{S_d})}|_{{S_{d_1}}|_{S_d}}}({({S_{d_2}}|_{S_d})}|_{{S_{d_1}}|_{S_d}})\\
=&\ 1-x-\frac{1}{dd_1(n+3-d-d_1-d_2x-d_3y)^n}\int_0^\frac{{n+3-d-d_1-d_2x-d_3y}}{d_2}dd_1(n+3-d-d_1-d_2x-d_3y-d_2t)^n\dt \\
=&\ 1-x-\frac{n+3-d-d_1-d_2x-d_3y}{d_2(n+1)}\\
=&\ \frac{d_2(n+1)-(n+3-d-d_1)}{d_2(n+1)}-\frac{nx}{n+1}+\frac{d_3y}{d_2(n+1)}\geq 0,
\end{align*}
\begin{align*}
&\ \beta_{{S_{d_1}}|_{S_d}, x{({S_{d_2}}|_{S_d})}|_{{S_{d_1}}|_{S_d}}+y{({S_{d_3}}|_{S_d})}|_{{S_{d_1}}|_{S_d}}}({({S_{d_3}}|_{S_d})}|_{{S_{d_1}}|_{S_d}})\\
=&\ 1-y-\frac{1}{dd_1(n+3-d-d_1-d_2x-d_3y)^n}\int_0^\frac{{n+3-d-d_1-d_2x-d_3y}}{d_3}dd_1(n+3-d-d_1-d_2x-d_3y-d_3t)^n\dt \\
=&\ 1-y-\frac{n+3-d-d_1-d_2x-d_3y}{d_3(n+1)}\\
=&\ \frac{d_3(n+1)-(n+3-d-d_1)}{d_3(n+1)}-\frac{ny}{n+1}+\frac{d_2x}{d_3(n+1)}\geq 0.
\end{align*}
We denote by $P$ the polytope generated by the following equations
\begin{equation*}
\begin{cases}
0\leq x\leq 1\\
0\leq y\leq 1\\
d_2x+d_3y\leq n+3-d-d_1\\
\frac{d_2(n+1)-(n+3-d-d_1)}{d_2(n+1)}-\frac{nx}{n+1}+\frac{d_3y}{d_2(n+1)}\geq 0\\
\frac{d_3(n+1)-(n+3-d-d_1)}{d_3(n+1)}-\frac{ny}{n+1}+\frac{d_2x}{d_3(n+1)}\geq 0.\\
\end{cases}
\end{equation*}
If $d_2+d_3<n+3-d-d_1$, then $P$ is generated by the extremal points 
$$(0,0), \ \ \left(\frac{(n+1)d_2-(n+3-d-d_1)}{nd_2}, 0\right),\ \ \left(0, \frac{(n+1)d_3-(n+3-d-d_1)}{nd_3}\right)$$
and
$$\left(\frac{nd_2+d_1+d_3-(n+3-d)}{(n-1)d_2}, \frac{nd_3+d_1+d_2-(n+3-d)}{(n-1)d_3}\right). $$
If $d_2+d_3\geq n+3-d-d_1$, then $P$ is generated by the extremal points
$$(0,0), \ \ \left(\frac{(n+1)d_2-(n+3-d-d_1)}{nd_2}, 0\right),\ \ \left(0, \frac{(n+1)d_3-(n+3-d-d_1)}{nd_3}\right)$$
and 
$$\left(1, \frac{n+3-d-d_1-d_2}{d_3}\right),\ \  \left(\frac{n+3-d-d_1-d_3}{d_2}, 1\right). $$
In either case, these points correspond to K-semistable log pairs by Theorem \ref{thm: SSS|S} and Remark \ref{rem: facts}, thus we exactly have 
$$\Kss\left({S_{d_1}}|_{S_d}, {({S_{d_2}}|_{S_d})}|_{{S_{d_1}}|_{S_d}}+{({S_{d_3}}|_{S_d})}|_{{S_{d_1}}|_{S_d}}\right)=P.$$
\end{example}

\
\
\

\section{$(\bP^n, S_{d_1}+S_{d_2}+S_{d_3})$}

In this section, we want to compute $\Kss(\bP^n, S_{d_1}+S_{d_2}+S_{d_3})$, where $S_{d_1}, S_{d_2}, S_{d_3}$ are general smooth hypersurfaces in $\bP^n$ with $d_1, d_2, d_3\leq n+1$ and $n\geq 2$.

\begin{theorem}\label{thm: PSSS}
Let $d_1, d_2, d_3$ be three positive integers satisfying $d_1+d_2+d_3<n+1$ for $n\geq 2$. 
Suppose all the Fano complete intersections are K-semistable.
%Suppose $S_{d_1}\cap S_{d_2}\cap S_{d_3}$ is K-semistable of dimension $n-3$. 
Then the log Fano pair
$$(\bP^n, aS_{d_1} +bS_{d_2}+cS_{d_3})$$
is K-semistable for
$$a=\frac{d_2+d_3+(n-1)d_1-(n+1)}{(n-2)d_1},$$
$$b=\frac{d_1+d_3+(n-1)d_2-(n+1)}{(n-2)d_2},$$
$$c=\frac{d_1+d_2+(n-1)d_3-(n+1)}{(n-2)d_3}. $$
\end{theorem}

\begin{proof}
By  Remark \ref{rem: kss degeneration},  $S_{d_1}$ induces a test configuration of $(\bP^n, aS_{d_1}+bS_{d_2}+cS_{d_3})$ such that the central fiber, denoted by $(X, aD_1+bD_2+cD_3)$, is the projective cone over $$(S_{d_1}, b{S_{d_2}}|_{S_{d_1}}+ c{S_{d_3}}|_{S_{d_1}})$$  
with respect to the polarization $\mO_{S_{d_1}}(d_1)$, and  $D_1\cong S_{d_1}$ is the infinite divisor. It suffices to show that $(X, aD_1+bD_2+cD_3)$ is a K-semistable log Fano pair.
We have the following computation:
$$\mO_{S_{d_1}}(d_1)\sim_\bQ -\frac{1}{r}(K_{S_{d_1}}+b{S_{d_2}}|_{S_{d_1}}+ c{S_{d_3}}|_{S_{d_1}})
\quad \text{and}  \quad
a=1-\frac{r}{n}$$
for 
$$r= \frac{n+1-d_1-bd_2-cd_3}{d_1}.$$ 
Note that $(S_{d_1}, b{S_{d_2}}|_{S_{d_1}}+ c{S_{d_3}}|_{S_{d_1}})$ 
is K-semistable by Theorem \ref{thm: SSS}. Applying  Lemma \ref{lem:cone stability}, we see that 
$(X, aD_1 +bD_2+cD_3)$
 is also K-semistable.
\end{proof}

\begin{example}\label{exa: PSSS}
Consider the log pair $(\bP^n, S_{d_1}+S_{d_2}+S_{d_3})$, where $S_{d_1}, S_{d_2}, S_{d_3}$ are three general smooth hypersurfaces of degrees $d_1, d_2, d_3$ in $\bP^{n}$ with $n\geq 2$ and $d_1, d_2, d_3\leq n+1$. 
\iffalse
Suppose the following two conditions:
\begin{enumerate}
\item $(\bP^n, S_{d_1}+S_{d_2}+S_{d_3})$ is log canonical,
\item $S_{d_1}\cap S_{d_2}\cap S_{d_3}$ is K-semistable of dimension $n-3$ if $d_1+d_2+d_3\leq n+1$. 
\end{enumerate}
\fi
Suppose all the Fano complete intersections are K-semistable.
%Suppose $S_{d_1}\cap S_{d_2}\cap S_{d_3}$ is K-semistable of dimension $n-3$ if $d_1+d_2+d_3\leq n+1$. 
We want to compute $\Kss(\bP^n, S_{d_1}+S_{d_2}+S_{d_3})$.

Suppose $(x,y,z)\in \Kss(\bP^n, S_{d_1}+S_{d_2}+S_{d_3})$, then $(\bP^n, xS_{d_1}+yS_{d_2}+zS_{d_3})$ is K-semistable. Applying beta criterion, we have
\begin{align*}
&\ \beta_{\bP^n, xS_{d_1}+yS_{d_2}+zS_{d_3}}(S_{d_1})\\
=&\ 1-x-\frac{1}{(n+1-d_1x-d_2y-d_3z)^n}\int_0^\frac{{n+1-d_1x-d_2y-d_3z}}{d_1}(n+1-d_1x-d_2y-d_3z-d_1t)^n\dt \\
=&\ 1-x-\frac{n+1-d_1x-d_2y-d_3z}{d_1(n+1)}\\
=&\ \frac{d_1-1}{d_1}-\frac{nx}{n+1}+\frac{d_2y}{d_1(n+1)}+\frac{d_3z}{d_1(n+1)}\geq 0,
\end{align*}
\begin{align*}
&\ \beta_{\bP^n, xS_{d_1}+yS_{d_2}+zS_{d_3}}(S_{d_2})\\
=&\ 1-y-\frac{1}{(n+1-d_1x-d_2y-d_3z)^n}\int_0^\frac{{n+1-d_1x-d_2y-d_3z}}{d_2}(n+1-d_1x-d_2y-d_3z-d_2t)^n\dt \\
=&\ 1-y-\frac{n+1-d_1x-d_2y-d_3z}{d_2(n+1)}\\
=&\ \frac{d_2-1}{d_2}-\frac{ny}{n+1}+\frac{d_1x}{d_2(n+1)}+\frac{d_3z}{d_2(n+1)}\geq 0,
\end{align*}
\begin{align*}
&\ \beta_{\bP^n, xS_{d_1}+yS_{d_2}+zS_{d_3}}(S_{d_3})\\
=&\ 1-z-\frac{1}{(n+1-d_1x-d_2y-d_3z)^n}\int_0^\frac{{n+1-d_1x-d_2y-d_3z}}{d_3}(n+1-d_1x-d_2y-d_3z-d_3t)^n\dt \\
=&\ 1-z-\frac{n+1-d_1x-d_2y-d_3z}{d_3(n+1)}\\
=&\ \frac{d_3-1}{d_3}-\frac{nz}{n+1}+\frac{d_1x}{d_3(n+1)}+\frac{d_2y}{d_3(n+1)}\geq 0.
\end{align*}

We denote by $P$ the polytope generated by the following equations
\begin{equation*}
\begin{cases}
0\leq x\leq 1\\
0\leq y\leq 1\\
0\leq z\leq 1\\
d_1x+d_2y+d_3z\leq n+1\\
\frac{d_1-1}{d_1}-\frac{nx}{n+1}+\frac{d_2y}{d_1(n+1)}+\frac{d_3z}{d_1(n+1)}\geq 0\\
\frac{d_2-1}{d_2}-\frac{ny}{n+1}+\frac{d_1x}{d_2(n+1)}+\frac{d_3z}{d_2(n+1)}\geq 0
\\
\frac{d_3-1}{d_3}-\frac{nz}{n+1}+\frac{d_1x}{d_3(n+1)}+\frac{d_2y}{d_3(n+1)}\geq 0.
\end{cases}
\end{equation*}
If $d_1+d_2+d_3<n+1$, then $P$ is generated by the extremal points 
$$(0,0,0), \ \ \left(\frac{(n+1)(d_1-1)}{nd_1}, 0,0\right),\ \ \left(0, \frac{(n+1)(d_2-1)}{nd_2}, 0\right), \ \ \left(0,0, \frac{(n+1)(d_3-1)}{nd_3}\right),$$

$$\left(0, \frac{d_3-1+n(d_2-1)}{(n-1)d_2}, \frac{d_2-1+n(d_3-1)}{(n-1)d_3}\right),$$

$$\left(\frac{d_3-1+n(d_1-1)}{(n-1)d_1}, 0,  \frac{d_1-1+n(d_3-1)}{(n-1)d_3}\right),$$

$$\left(\frac{d_2-1+n(d_1-1)}{(n-1)d_1}, \frac{d_1-1+n(d_2-1)}{(n-1)d_2}, 0\right),
$$

$$\left(\frac{d_2+d_3+(n-1)d_1-(n+1)}{(n-2)d_1}, \frac{d_1+d_3+(n-1)d_2-(n+1)}{(n-2)d_2}, \frac{d_1+d_2+(n-1)d_3-(n+1)}{(n-2)d_3}\right). $$
If $d_1+d_2+d_3\geq n+1$, then $P$ is generated by the extremal points
$$(0,0,0), \ \ \left(\frac{(n+1)(d_1-1)}{nd_1}, 0,0\right),\ \ \left(0, \frac{(n+1)(d_2-1)}{nd_2}, 0\right), \ \ \left(0,0, \frac{(n+1)(d_3-1)}{nd_3}\right),$$
and $A_1 ,B_1, A_2 ,B_2, A_3 ,B_3$, where
\begin{equation*}
A_1, B_1=
\begin{cases}
(0, 1, \frac{n+1-d_2}{d_3}),\ (0,\frac{n+1-d_3}{d_2}, 1) & \text{if \ $d_2+d_3\geq n+1$}\\
\ \\
\left(0, \frac{d_3-1+n(d_2-1)}{(n-1)d_2}, \frac{d_2-1+n(d_3-1)}{(n-1)d_3}\right),\ \left(\frac{n+1-(d_2+d_3)}{d_1}, 1,1\right) & \text{if \ $d_2+d_3\leq n+1$}\\
\end{cases}
\end{equation*}

\begin{equation*}
A_2, B_2=
\begin{cases}
(1, 0, \frac{n+1-d_1}{d_3}),\ (\frac{n+1-d_3}{d_1}, 0, 1) & \text{if \ $d_1+d_3\geq n+1$}\\
\ \\
\left(\frac{d_3-1+n(d_1-1)}{(n-1)d_1}, 0,  \frac{d_1-1+n(d_3-1)}{(n-1)d_3}\right),\ \left(1, \frac{n+1-(d_1+d_3)}{d_2}, 1\right) & \text{if \ $d_1+d_3\leq n+1$}\\
\end{cases}
\end{equation*}

\begin{equation*}
A_3, B_3=
\begin{cases}
(1, \frac{n+1-d_1}{d_2}, 0),\ (\frac{n+1-d_2}{d_1}, 1, 0) & \text{if \ $d_1+d_2\geq n+1$}\\
\ \\
\left(\frac{d_2-1+n(d_1-1)}{(n-1)d_1}, \frac{d_1-1+n(d_2-1)}{(n-1)d_2}, 0\right),\ \left(1,1,\frac{n+1-(d_1+d_2)}{d_3}\right) & \text{if \ $d_1+d_2\leq n+1$}.\\
\end{cases}
\end{equation*}
In either case, these points correspond to K-semistable log pairs by Theorem \ref{thm: PSSS}, Remark \ref{rem: facts}, and Theorem \ref{thm: PSS}, thus we exactly have 
$$\Kss(\bP^n, S_{d_1}+S_{d_2}+S_{d_3})=P.$$
\end{example}

\begin{example}
In Example \ref{exa: PSSS}, take $n\geq 3, d_1=2, d_2=d_3=1$. We denote $S_{d_1}, S_{d_2}, S_{d_3}$ by $Q, H_1, H_2$, then $\Kss(\bP^n, Q+H_1+H_2)$ is generated by the extremal points
$$(0,0,0),\ (\frac{n+1}{2n}, 0,0),\ \left(\frac{n}{2(n-1)}, 0, \frac{1}{n-1}\right),\ \left(\frac{n}{2(n-1)}, \frac{1}{n-1},0\right),\ \left(\frac{n-1}{2(n-2)},\frac{1}{n-2}, \frac{1}{n-2}\right). $$ 
\end{example}

\begin{example}
In Example \ref{exa: PSSS}, take $n= 3, d_1=d_2=2, d_3=1$. We denote $S_{d_1}, S_{d_2}, S_{d_3}$ by $Q, Q', H$, then $\Kss(\bP^3, Q+Q'+H)$ is generated by the extremal points
$$(0,0,0),\ (\frac{2}{3}, 0,0),\ (0, \frac{2}{3},0),\ (0, \frac{3}{4}, \frac{1}{2}),\ (\frac{1}{2},1,1),\  (\frac{3}{4}, 0, \frac{1}{2}),\ (1, \frac{1}{2},1),\ (1,1,0). $$ 
\end{example}

\begin{example}
In Example \ref{exa: PSSS}, take $n=3, d_1=d_2=d_3=2$. We denote $S_{d_1}, S_{d_2}, S_{d_3}$ by $Q, Q', Q''$, then $\Kss(\bP^3, Q+Q'+Q'')$ is generated by the extremal points
$$(0,0,0),\ (\frac{2}{3}, 0,0),\ (0, \frac{2}{3}, 0),\ (0,0, \frac{2}{3}),\ (0,1,1),\ (1,0,1),\ (1,1,0). $$ 
\end{example}

\
\
\

\section{$(S_d, {S_{d_1}}|_{S_d}+{S_{d_2}}|_{S_d}+{S_{d_3}}|_{S_d})$}

In this section, we want to compute $\Kss(S_d, {S_{d_1}}|_{S_d}+{S_{d_2}}|_{S_d}+{S_{d_3}}|_{S_d})$, where $S_d, S_{d_1}, S_{d_2}, S_{d_3}$ are general smooth hypersurfaces in $\bP^{n+1}$ with $d_1, d_2, d_3\leq n+2-d$ and $n\geq 2$.

\begin{theorem}\label{thm: SSSS}
Let $d, d_1, d_2, d_3$ be four positive integers satisfying $d_1+d_2+d_3<n+2-d$ for $n\geq 2$. 
Suppose all the Fano complete intersections are K-semistable.
%Suppose $S_d\cap S_{d_1}\cap S_{d_2}\cap S_{d_3}$ is K-semistable of dimension $n-3$. 
Then the log Fano pair
$$(S_d, a{S_{d_1}}|_{S_d} +b{S_{d_2}}|_{S_d}+c{S_{d_3}}|_{S_d})$$
is K-semistable for
$$a=\frac{d_2+d_3+(n-1)d_1-(n+2-d)}{(n-2)d_1},$$
$$b=\frac{d_1+d_3+(n-1)d_2-(n+2-d)}{(n-2)d_2},$$
$$c=\frac{d_1+d_2+(n-1)d_3-(n+2-d)}{(n-2)d_3}. $$
\end{theorem}

\begin{proof}
By  Remark \ref{rem: kss degeneration},  $S_{d_1}\cap S_d$ induces a test configuration of 
$$(S_d, a{S_{d_1}}|_{S_d} +b{S_{d_2}}|_{S_d}+c{S_{d_3}}|_{S_d})$$ 
such that the central fiber, denoted by $(X, aD_1+bD_2+cD_3)$, is the projective cone over 
$$(S_{d}\cap S_{d_1}, bS_d\cap S_{d_1}\cap S_{d_2}+cS_d\cap S_{d_1}\cap S_{d_3})= \left({S_{d_1}}|_{S_d}, b({S_{d_2}}|_{S_d})|_{{S_{d_1}}|_{S_d}}+ c({S_{d_3}}|_{S_d})|_{{S_{d_1}}|_{S_d}}\right)$$ with respect to the polarization $\mO_{S_d\cap S_{d_1}}(d_1)$, and  $D_1\cong S_d\cap S_{d_1}$ is the infinite divisor. 
It suffices to show that 
$(X, aD_1 +bD_2+cD_3)$
is a K-semistable log Fano pair.   
We have the following computation:
$$\mO_{{S_{d_1}}|_{S_d}}(d_1)\sim_\bQ -\frac{1}{r}\left(K_{{S_{d_1}}|_{S_d}}+b({S_{d_2}}|_{S_d})|_{{S_{d_1}}|_{S_d}}+ c({S_{d_3}}|_{S_d})|_{{S_{d_1}}|_{S_d}}\right)$$
and  
$$a=1-\frac{r}{n}$$
for 
$$r= \frac{n+2-d-d_1-bd_2-cd_3}{d_1}.$$ 
Note that $({S_{d_1}}|_{S_d}, b({S_{d_2}}|_{S_d})|_{{S_{d_1}}|_{S_d}}+ c({S_{d_3}}|_{S_d})|_{{S_{d_1}}|_{S_d}})$ 
is K-semistable by Theorem \ref{thm: SSS|S}.  Applying  Lemma \ref{lem:cone stability}, we see that 
$(X, aD_1 +bD_2+cD_3)$
 is also K-semistable.
\end{proof}

\begin{example}\label{exa: SSSS}
Consider the log pair $(S_d, {S_{d_1}}|_{S_d}+{S_{d_2}}|_{S_d}+{S_{d_3}}|_{S_d})$, where $S_d, S_{d_1}, S_{d_2}, S_{d_3}$ are four general smooth hypersurfaces in $\bP^{n+1}$ of degrees $d, d_1, d_2, d_3$ with $n\geq 2$ and $d_1, d_2, d_3\leq n+2-d$. 
\iffalse
Suppose the following two conditions:
\begin{enumerate}
\item $(S_d, {S_{d_1}}|_{S_d}+{S_{d_2}}|_{S_d}+{S_{d_3}}|_{S_d})$ is log canonical,
\item $S_d\cap S_{d_1}\cap S_{d_2}\cap S_{d_3}$ is K-semistable of dimension $n-3$ if $d_1+d_2+d_3\leq n+2-d$. 
\end{enumerate}
\fi
Suppose all the Fano complete intersections are K-semistable.
%Suppose $S_d\cap S_{d_1}\cap S_{d_2}\cap S_{d_3}$ is K-semistable of dimension $n-3$ if $d_1+d_2+d_3\leq n+2-d$.
We want to compute $\Kss(S_d, {S_{d_1}}|_{S_d}+{S_{d_2}}|_{S_d}+{S_{d_3}}|_{S_d})$.

Suppose $(x,y,z)\in \Kss(S_d, {S_{d_1}}|_{S_d}+{S_{d_2}}|_{S_d}+{S_{d_3}}|_{S_d})$, then $(S_d, x{S_{d_1}}|_{S_d}+y{S_{d_2}}|_{S_d}+z{S_{d_3}}|_{S_d})$ is K-semistable. Applying beta criterion, we have
\begin{align*}
&\ \beta_{S_d, x{S_{d_1}}|_{S_d}+y{S_{d_2}}|_{S_d}+z{S_{d_3}}|_{S_d}}({S_{d_1}}|_{S_d})\\
=&\ 1-x-\frac{1}{d(n+2-d-d_1x-d_2y-d_3z)^n}\int_0^\frac{{n+2-d-d_1x-d_2y-d_3z}}{d_1}d(n+2-d-d_1x-d_2y-d_3z-d_1t)^n\dt \\
=&\ 1-x-\frac{n+2-d-d_1x-d_2y-d_3z}{d_1(n+1)}\\
=&\ \frac{d_1(n+1)-(n+2-d)}{d_1(n+1)}-\frac{nx}{n+1}+\frac{d_2y}{d_1(n+1)}+\frac{d_3z}{d_1(n+1)}\geq 0,
\end{align*}
\begin{align*}
&\ \beta_{S_d, x{S_{d_1}}|_{S_d}+y{S_{d_2}}|_{S_d}+z{S_{d_3}}|_{S_d}}({S_{d_2}}|_{S_d})\\
=&\ 1-y-\frac{1}{d(n+2-d-d_1x-d_2y-d_3z)^n}\int_0^\frac{{n+2-d-d_1x-d_2y-d_3z}}{d_2}d(n+2-d-d_1x-d_2y-d_3z-d_2t)^n\dt \\
=&\ 1-y-\frac{n+2-d-d_1x-d_2y-d_3z}{d_2(n+1)}\\
=&\ \frac{d_2(n+1)-(n+2-d)}{d_2(n+1)}-\frac{ny}{n+1}+\frac{d_1x}{d_2(n+1)}+\frac{d_3z}{d_2(n+1)}\geq 0,
\end{align*}
\begin{align*}
&\ \beta_{S_d, x{S_{d_1}}|_{S_d}+y{S_{d_2}}|_{S_d}+z{S_{d_3}}|_{S_d}}({S_{d_3}}|_{S_d})\\
=&\ 1-z-\frac{1}{d(n+2-d-d_1x-d_2y-d_3z)^n}\int_0^\frac{{n+2-d-d_1x-d_2y-d_3z}}{d_3}d(n+2-d-d_1x-d_2y-d_3z-d_3t)^n\dt \\
=&\ 1-z-\frac{n+2-d-d_1x-d_2y-d_3z}{d_3(n+1)}\\
=&\ \frac{d_3(n+1)-(n+2-d)}{d_3(n+1)}-\frac{nz}{n+1}+\frac{d_1x}{d_3(n+1)}+\frac{d_2y}{d_3(n+1)}\geq 0.
\end{align*}

We denote by $P$ the polytope generated by the following equations
\begin{equation*}
\begin{cases}
0\leq x\leq 1\\
0\leq y\leq 1\\
0\leq z\leq 1\\
d_1x+d_2y+d_3z\leq n+2-d\\
\frac{d_1(n+1)-(n+2-d)}{d_1(n+1)}-\frac{nx}{n+1}+\frac{d_2y}{d_1(n+1)}+\frac{d_3z}{d_1(n+1)}\geq 0\\
\frac{d_2(n+1)-(n+2-d)}{d_2(n+1)}-\frac{ny}{n+1}+\frac{d_1x}{d_2(n+1)}+\frac{d_3z}{d_2(n+1)}\geq 0\\
\frac{d_3(n+1)-(n+2-d)}{d_3(n+1)}-\frac{nz}{n+1}+\frac{d_1x}{d_3(n+1)}+\frac{d_2y}{d_3(n+1)}\geq 0.
\end{cases}
\end{equation*}
If $d_1+d_2+d_3<n+2-d$, then $P$ is generated by the extremal points 
$$(0,0,0), \  \left(\frac{(n+1)d_1-(n+2-d)}{nd_1}, 0,0\right),\  \left(0, \frac{(n+1)d_2-(n+2-d)}{nd_2}, 0\right), \  \left(0,0, \frac{(n+1)d_3-(n+2-d)}{nd_3}\right),$$

$$\left(0, \frac{d_3+nd_2-(n+2-d)}{(n-1)d_2}, \frac{d_2+nd_3-(n+2-d)}{(n-1)d_3}\right),$$
  
$$\left(\frac{d_3+nd_1-(n+2-d)}{(n-1)d_1}, 0,  \frac{d_1+nd_3-(n+2-d)}{(n-1)d_3}\right), $$
 
$$\left(\frac{d_2+nd_1-(n+2-d)}{(n-1)d_1}, \frac{d_1+nd_2-(n+2-d)}{(n-1)d_2}, 0\right), $$

$$\left(\frac{d_2+d_3+(n-1)d_1-(n+2-d)}{(n-2)d_1}, \frac{d_1+d_3+(n-1)d_2-(n+2-d)}{(n-2)d_2}, \frac{d_1+d_2+(n-1)d_3-(n+2-d)}{(n-2)d_3}\right). $$
If $d_1+d_2+d_3\geq n+2-d$, then $P$ is generated by the extremal points
$$(0,0,0), \  \left(\frac{(n+1)d_1-(n+2-d)}{nd_1}, 0,0\right),\  \left(0, \frac{(n+1)d_2-(n+2-d)}{nd_2}, 0\right), \  \left(0,0, \frac{(n+1)d_3-(n+2-d)}{nd_3}\right),$$
and $A_1 ,B_1, A_2 ,B_2, A_3,B_3$, where
\begin{equation*}
A_1, B_1=
\begin{cases}
(0, 1, \frac{n+2-d-d_2}{d_3}),\ (0,\frac{n+2-d-d_3}{d_2}, 1) & \text{if \ $d_2+d_3\geq n+2-d$}\\
\ \\
\left(0, \frac{d_3+nd_2-(n+2-d)}{(n-1)d_2}, \frac{d_2+nd_3-(n+2-d)}{(n-1)d_3}\right),\ \left(\frac{n+2-d-(d_2+d_3)}{d_1},1,1\right) & \text{if \ $d_2+d_3\leq n+2-d$}\\
\end{cases}
\end{equation*}

\begin{equation*}
A_2, B_2=
\begin{cases}
(1, 0, \frac{n+2-d-d_1}{d_3}),\ (\frac{n+2-d-d_3}{d_1}, 0, 1) & \text{if \ $d_1+d_3\geq n+2-d$}\\
\ \\
\left(\frac{d_3+nd_1-(n+2-d)}{(n-1)d_1}, 0,  \frac{d_1+nd_3-(n+2-d)}{(n-1)d_3}\right),\ \left(1,\frac{n+2-d-(d_1+d_3)}{d_2},1\right) & \text{if \ $d_1+d_3\leq n+2-d$}\\
\end{cases}
\end{equation*}

\begin{equation*}
A_3, B_3=
\begin{cases}
(1, \frac{n+2-d-d_1}{d_2}, 0),\ (\frac{n+2-d-d_2}{d_1}, 1, 0) & \text{if \ $d_1+d_2\geq n+2-d$}\\
\ \\
\left(\frac{d_2+nd_1-(n+2-d)}{(n-1)d_1}, \frac{d_1+nd_2-(n+2-d)}{(n-1)d_2}, 0\right),\ \left(1,1, \frac{n+2-d-(d_1+d_2)}{d_3}\right) & \text{if \ $d_1+d_2\leq n+2-d$}.\\
\end{cases}
\end{equation*}
In either case, these points correspond to K-semistable log pairs by Theorem \ref{thm: SSSS}, Remark \ref{rem: facts}, and Theorem \ref{thm: SSS},  thus we exactly have 
$$\Kss(S_d, {S_{d_1}}|_{S_d}+{S_{d_2}}|_{S_d}+{S_{d_3}}|_{S_d})=P.$$
\end{example}

\begin{example}
In Example \ref{exa: SSSS}, take $n=3, d_1=d_2=d_3=1$ and denote $S_d, S_{d_1}, S_{d_2}, S_{d_3}$ by $Q, H_1, H_2, H_3$. Then $\Kss(Q, {H_1}|_Q+{H_2}|_Q+{H_3}|_Q)$ is generated by
$$(0,0,0),\ (\frac{1}{3},0,0),\ (0,\frac{1}{3},0),\ (0,0,\frac{1}{3}),\ (0,\frac{1}{2},\frac{1}{2}),\ (\frac{1}{2},0,\frac{1}{2}),\ (\frac{1}{2},\frac{1}{2},0),\ (1,1,1). $$
\end{example}

\
\
\

\section{Final remark}

From the computation in the previous sections, we see that the same method applies to the case where there are arbitrary number of components. In this section, we prove the following general result.

\begin{theorem}\label{thm: k components}
Consider the log pair $(\bP^n, S_{d_1}+S_{d_2}+...+S_{d_k})$, where $S_{d_j}, j=1,...,k,$ are general smooth hypersurfaces in $\bP^n$ of degrees $d_j$ with $n\geq 2$ and $d_j\leq n+1$. 
\iffalse
Suppose the following two conditions:
\begin{enumerate}
\item $(\bP^n, S_{d_1}+...+S_{d_k})$ is log canonical;
\item $S_{d_1}\cap...\cap S_{d_k}$ is K-semistable of dimension $n-k$ if $\sum_{j=1}^kd_j\leq n+1$. 
\end{enumerate}
\fi
Suppose all the Fano complete intersections are K-semistable.
Then $\Kss(\bP^n, S_{d_1}+...+S_{d_k})$ is a polytope generated by the following equations
\begin{equation*}
\begin{cases}
0\leq x_i\leq 1, \ \ \  1\leq i\leq k\\
\ \\
\beta_{\bP^n, \sum_{j=1}^kx_jS_{d_j}}(S_{d_i})\geq 0, \ \ \  1\leq i\leq k\\
\ \\
\sum_{j=1}^k x_jd_j\leq n+1.
\end{cases}
\end{equation*}
\end{theorem}

\begin{proof}
We denote by $P$ the polytope given by the above equations, and aim to show that all the extremal points of $P$ correspond to K-semistable log pairs. We divide the proof into two steps.

\

\textit{Step 1}.
We first assume $\sum_{j=1}^k d_j<n+1$, then the equations 
$$\beta_{\bP^n, \sum_{j=1}^kx_jS_{d_j}}(S_{d_i})= 0, \ \ \ 1\leq i\leq k$$ 
admit a unique solution, denoted by $(a_1,...,a_k)$, where
$$a_j=\frac{\sum_{i\ne j}d_i+(n-k+2)d_j-(n+1)}{(n-k+1)d_j}. $$
Clearly $(a_1,...,a_k)$ is an extremal point of $P$. We first show that $(\bP^n, \sum_{j=1}^k a_jS_{d_j})$ is a K-semistable log Fano pair.

By  Remark \ref{rem: kss degeneration},  $S_{d_1}$ induces a test configuration of $(\bP^n, \sum_{j=1}^ka_jS_{d_j})$ such that the central fiber, denoted by $(X, \sum_{j=1}^ka_jD_j)$, is the projective cone over $(S_{d_1}, \sum_{j=2}^k a_jS_{d_1}\cap S_{d_j})$ with respect to the polarization $\mO_{S_{d_1}}(d_1)$, and  $D_1\cong S_{d_1}$ is the infinite divisor. It suffices to show that $(X, \sum_{j=1}^ka_jD_j)$ is a K-semistable log Fano pair.
We have the following computation:
$$\mO_{S_{d_1}}(d_1)\sim_\bQ -\frac{1}{r}(K_{S_{d_1}}+\sum_{j=2}^ka_j{S_{d_j}}|_{S_{d_1}})
\quad \text{and}  \quad
a_1=1-\frac{r}{n}$$
for 
$$r= \frac{n+1-d_1-\sum_{j=2}^ka_jd_j}{d_1}.$$ 
By Lemma \ref{lem:cone stability}, it suffices to show that 
$$(S_{d_1}, \sum_{j=2}^ka_j{S_{d_j}}|_{S_{d_1}})$$ 
is K-semistable.  Applying the same degeneration approach as above, it suffices to show that
$$\left({S_{d_2}}|_{S_{d_1}}, \sum_{j=3}^k a_j\left({S_{d_j}}|_{S_{d_1}}\right)|_{{S_{d_2}}|_{S_{d_1}}}\right)=(S_{d_2}\cap S_{d_1}, \sum_{j=3}^k a_jS_{d_j}\cap S_{d_2}\cap S_{d_1}) $$
is K-semistable. Step by step, one could finally reduce it to show that
$$(S_{d_{k-1}}\cap ... \cap S_{d_1}, a_k S_{d_k}\cap S_{d_{k-1}}\cap ... \cap S_{d_1}) $$
is K-semistable, which could be derived by applying Lemma \ref{lem: zz general}.

Since we assume $\sum_{j=1}^k d_j<n+1$, the other extremal points must appear on the hyperplanes $\{x_j=0\}_{j=1}^k$. For each such hyperplane, one could apply Theorem \ref{thm: k components} where there are $k-1$ components in the boundary (we have proved Theorem \ref{thm: k components} for $k=1,2,3$). Thus we see that all the extremal points of $P$ correspond to K-semistable log pairs.

\

\textit{Step 2}. In this step, we assume $\sum_{j=1}^k d_j\geq n+1$. For this case, there could be some extremal points appearing on the hyperplane $\sum_{j=1}^k d_j x_j=n+1$. If $(b_1,...,b_k)$ is such an extremal point of $P$, then it is clear that $(\bP^n, \sum_{j=1}^k b_j S_{d_j})$ is an lc log Calabi-Yau pair, which is K-semistable. For other type of extremal points of $P$, the same as before, they appear on the hyperplanes $\{x_j=0\}_{j=1}^k$. We are done by applying Theorem \ref{thm: k components} where there are $k-1$ components in the boundary. Thus all the extremal points of $P$ correspond to K-semistable log pairs.
\end{proof}

Similarly, we have the following result.

\begin{theorem}\label{thm: k components'}
Consider the log pair $(S_{d}, {S_{d_1}}|_{S_d}+{S_{d_2}}|_{S_d}+...+{S_{d_k}}|_{S_d})$, where $S_d, S_{d_j}, j=1,...,k,$ are general smooth hypersurfaces in $\bP^{n+1}$ of degrees $d, d_j$ with $n\geq 2$ and $d_j\leq n+2-d$. 
\iffalse
Suppose the following two conditions:
\begin{enumerate}
\item $(S_d, {S_{d_1}}|_{S_d}+...+{S_{d_k}}|_{S_d})$ is log canonical;
\item $S_d\cap S_{d_1}\cap...\cap S_{d_k}$ is K-semistable of dimension $n-k$ if $\sum_{j=1}^kd_j\leq n+2-d$. 
\end{enumerate}
\fi
Suppose all the Fano complete intersections are K-semistable.
%Suppose $S_d\cap S_{d_1}\cap...\cap S_{d_k}$ is K-semistable of dimension $n-k$ if $\sum_{j=1}^kd_j\leq n+2-d$. 
Then $\Kss(S_{d}, {S_{d_1}}|_{S_d}+{S_{d_2}}|_{S_d}+...+{S_{d_k}}|_{S_d})$ is a polytope generated by the following equations
\begin{equation*}
\begin{cases}
0\leq x_i\leq 1, \ \ \  1\leq i\leq k\\
\ \\
\beta_{S_d, \sum_{j=1}^k x_j{S_{d_j}}|_{S_d}}({S_{d_i}}|_{S_d})\geq 0, \ \ \  1\leq i\leq k\\
\ \\
\sum_{j=1}^k x_jd_j\leq n+2-d.
\end{cases}
\end{equation*}
\end{theorem}

We end this article by formulating the following conjecture, which we will study in a future work.

\begin{conjecture}
Let $(X, \sum_{j=1}^k D_k)$ be a log pair taken from $\mE:=\mE(d,k, v, I)$ (see Section \ref{sec: preliminaries}), then
\begin{enumerate}
\item $\Kss(X, \sum_{j=1}^kD_j)$ is a polytope in the simplex $\overline{\Delta^k}$ (see Section \ref{sec: preliminaries});
\item all the extremal points of $\Kss(X, \sum_{j=1}^kD_j)$ are rational, which means that their coordinates are rational numbers;
\item all the extremal points of $\Kss(X, \sum_{j=1}^kD_j)$ depend only on $d,k,v,I$.
\end{enumerate}
\end{conjecture}

\begin{remark}
All the examples we study in this article satisfy this conjecture. We also mention that the conjecture for $k=1$ was completely answered in \cite{Zhou23}.
\end{remark}

\bibliography{reference.bib}
\end{document}